\documentclass[a4paper,12pt]{article} 
\usepackage[english]{babel}
\usepackage{graphicx,latexsym,amsmath,amsthm,amssymb,color,indentfirst,xfrac}
 \usepackage[usenames,dvipsnames]{pstricks}
\usepackage{mathrsfs} 
 \usepackage{epsfig}
\usepackage{pst-grad} 
 \usepackage{pst-plot} 

\newtheorem{satz}{thm1}

\newtheorem{sat}{thm2}

\newtheorem{thm}[sat]{Theorem}
\newtheorem{lem}[satz]{Lemma}
\newtheorem{cor}[satz]{Corollary}
\newtheorem{prop}[satz]{Proposition}
\newtheorem{ex}[satz]{Example}
\theoremstyle{definition}
\newtheorem{defn}[satz]{Definition}

\newcommand{\be}{\begin{equation}}           
\newcommand{\ee}{\end{equation}}

\newcommand{\ba}{\begin{align}}                
\newcommand{\ea}{\end{align}}

\newcommand{\bal}{\begin{align*}}              
\newcommand{\eal}{\end{align*}}

\newcommand{\bxx}{\begin{ex}}

\newcommand{\exx}{\end{ex}}

\newcommand{\txsit}

\newenvironment{pr}
{\begin{trivlist}
\item[\hskip\labelsep{\bf Proof.}]}                     
{$\hfill\Box$\end{trivlist}}

\title{Every planar graph without $i$-cycles adjacent simultaneously to $j$-cycles and $k$-cycles is DP-$4$-colorable when $\{i,j,k\}=\{3,4,5\}$}
\author {Pongpat Sittitrai\\ 
{\small\em Department of Mathematics, Faculty of Science, Khon Kaen University, 40002, Thailand }\\  
{\small\em E-mail address: pongpat\_s@kkumail.com} 
\and Kittikorn Nakprasit \footnote{Corresponding Author} \\ 
{\small\em Department of Mathematics, Faculty of Science, Khon Kaen University, 40002, Thailand }\\
{\small\em E-mail address: kitnak@hotmail.com}}
\date{}

\begin{document}

\maketitle

\begin{center}{\bf Abstract}\end{center}
\indent\indent

DP-coloring is a generalization of a list coloring in simple graphs. 
Many results in list coloring can be generalized in those of DP-coloring. 
Kim and Ozeki showed that planar graphs without $k$-cycles 
where $k=3,4,5,$ or $6$ are DP-$4$-colorable. 
Recently, Kim and Yu extended the result on $3$- and $4$-cycles by showing 
that planar graphs without triangles adjacent to $4$-cycles are DP-$4$-colorable. 
Xu and Wu showed that planar graphs without $5$-cycles adjacent 
simultaneously to $3$-cycles and $4$-cycles are $4$-choosable. 
In this paper, we extend the result on $5$-cycles  
and triangles adjacent to $4$-cycles  
by showing that planar graphs without $i$-cycles adjacent simultaneously 
to $j$-cycles and $k$-cycles are DP-$4$-colorable when $\{i,j,k\}=\{3,4,5\}.$ 
This also generalizes the result of Xu and Wu. 

\section{Introduction}
Every graph in this paper is finite, simple, and undirected. 
Embedding a graph $G$ in the plane, we let $V(G), E(G),$ and $F(G)$ 
denote the vertex set, edge set, and face set of $G.$ 
For $U \subseteq V(G),$ we let $G[U]$ denote 
the subgraph of $G$ induced by $U.$ For $X, Y \subseteq V(G)$ 
where $X$ and $Y$ are disjoint, we let $E_G(X,Y)$ be the set 
of all edges in $G$ with one endpoint in $X$ and the other in $Y.$ 

The concept of choosability was introduced by Vizing in 1976 \cite{Vizing} 
and by Erd\H os, Rubin, and Taylor in 1979 \cite{Erdos}, independently. 
A $k$-\emph{list assignment} $L$ of a graph $G$ assigns a list $L(v)$ 
(a set of colors) with $|L(v)|= k$ to each vertex $v.$ 
A graph $G$ is $L$-colorable if there is a proper coloring $f$ where $f(v)\in L(v).$   
If $G$ is $L$-colorable for every $k$-assignment $L,$ then we say  $G$ 
is $k$-\emph{choosable}. 

Dvo\v{r}\'{a}k and Postle \cite{DP} introduced a generalization 
of list coloring in which they called a \emph{correspondence coloring}. 
But following Bernshteyn, Kostochka, and Pron \cite{BKP},  
we call it a \emph{DP-coloring}.   


\begin{defn}\label{H}
Let $L$ be an assignment of a graph $G.$ 
We call $H$ a \emph{cover} of $G$ 
if it satisfies all the followings:\\ 
(i) The vertex set of $H$ is $\bigcup_{u \in V(G)}(\{u\}\times L(u)) =
\{(u,c): u \in V(G), c \in L(u) \};$\\
(ii) $H[u\times L(u)]$ is a complete graph for every $u \in V(G);$\\
(iii) For each $uv \in E(G),$ 
the set $E_H(\{u\}\times L(u), \{v\}\times L(v))$ is a matching (maybe empty). 
(iv) If $uv \notin E(G),$ then no edges of $H$ connect  
$\{u\}\times L(u)$ and  $\{v\}\times L(v).$   
\end{defn}

\begin{defn}\label{DP} 
An $(H,L)$-coloring of $G$ is an independent set in 
a cover $H$ of $G$ with size $|V(G)|.$ 
We say that a graph is \emph{DP-$k$-colorable} if $G$ has 
an $(H,L)$-coloring for every $k$-assignment $L$ and every cover $H$ of $(G.$   
The \emph{DP-chromatic number} of $G,$ denoted by $\chi_{DP}(G),$ 
is the minimum number $k$ such that $G$ is DP-$k$-colorable. 
\end{defn}

If we define edges on $H$ to match exactly the same colors in $L(u)$ and $L(v)$ 
for each $uv \in E(H),$
then $G$ has an $(H,L)$-coloring if and only 
if $G$ is $L$-colorable.  
Thus DP-coloring is a generalization of list coloring. 
This also implies that  $\chi_{DP}(G) \geq \chi_l(G).$ 
In fact, the difference of these two chromatic numbers can be arbitrarily large. 
For graphs with average degree $d,$ Bernshteyn \cite{Bern} 
showed that   $\chi_{DP}(G) = \Omega(d /\log d),$   
while Alon \cite{Alon} showed that   $\chi_l(G) = \Omega(\log d).$ 

Dvo\v{r}\'{a}k  and Postle \cite{DP} showed that 
$\chi_{DP}(G) \leq 5$ for every planar graph $G.$  
This extends a seminal result by Thomassen \cite{Tho} on list colorings. 
On the other hand, Voigt \cite{Vo1}  gave an example of a planar graph 
which is not $4$-choosable (thus not DP-$4$-colorable). 
It is of interest to obtain sufficient conditions 
for planar graphs to be DP-$4$-colorable. 
Kim and Ozeki \cite{KimO} showed that planar graphs without $k$-cycles 
are DP-$4$-colorable for each $k =3,4,5,6.$  
Kim and Yu  \cite{KimY} extended the result on $3$- and $4$-cycles 
by showing that planar graphs without triangles adjacent to $4$-cycles  
are DP-$4$-colorable. 
In this paper, we extend the result on $5$-cycles and triangles adjacent to $4$-cycles  
by showing that planar graphs without $i$-cycles adjacent simultaneously to $j$-cycles 
and $k$-cycles are DP-$4$-colorable when $\{i,j,k\}=\{3,4,5\}.$

\begin{thm}\label{main2}
Every planar graph without $i$-cycles adjacent simultaneously 
to $j$-cycles and $k$-cycles is DP-$4$-colorable when $\{i,j,k\}=\{3,4,5\}$.
\end{thm}

Theorem~\ref{main2} generalizes the result of Xu and Wu \cite{Xu} as follows. 
\begin{cor}
Every planar graphs without $5$-cycles adjacent simultaneously 
to $3$-cycles and $4$-cycles is $4$-choosable. 
\end{cor}

\section{Structure Obtained from Condition on Cycles}
First, we introduce some notations and definitions. 
A $k$-vertex ($k^+$-vertex, $k^-$-vertex, respectively) is 
a vertex of degree $k$ (at least $k,$ at most $k,$ respectively. 
The same notations are applied to faces. 
A \emph{$(d_1,d_2,\dots,d_k)$-face} $f$ is a face of degree $k$ 
where all vertices on $f$ have degree $d_1,d_2,\dots,d_k$ in any arbitrary order. 
A \emph{$(d_1,d_2,\dots,d_k)$-vertex} $v$ is a vertex of degree $k$ 
where all faces incident to $v$ have degree $d_1,d_2,\dots,d_k$ in any arbitrary order.  
A graph $C(m,n)$ is a plane graph obtained from an $(m+n-2)$-cycle with one chord 
such that internal faces have length $m$ and $n.$ 
A graph $C(l,m,n)$ is a plane graph obtained from an $(l+m+n-4)$-cycle with two chords 
such that internal faces have length $l,m,$ and $n$ 
where a middle  face has length $m.$  
Note that each $C(l,m,n)$ is not necessary unique. 
For example, there are two non isomorphic graphs that are $C(3,4,3).$

Let $G$ be a graph without $i$-cycles adjacent simultaneously to $j$-cycles and $k$-cycles  where $\{i,j,k\}=\{3,4,5\}$ 
The following property is straightforward. 
\begin{prop}\label{prop1} 
	$G$  does not contain $C(3,4), C(3,3,3),$ and $C(3,3,5).$ 
\end{prop}

Proposition \ref{prop1} yields the following immediate consequence.  

\begin{prop} \label{prop11} 
	If $v$ be an $n$-vertex in $G,$ then $v$ in $G$ is incident to at most $n-2$ $3$-faces.
\end{prop} 

\section {Structure of Minimal Non DP-$4$-colorable Graphs}

\begin{defn}\label{residual} 
	Let $H$ be a cover of $G$ with list assignment $L.$ 
	Let $G'= G-F$ where $F$ is an induced subgraph of $G.$ 
	A list assignment $L'$ is a \emph{restriction of $L$} on $G'$ 
	if $L'(u) = L(u)$ for each vertex in $G'.$ 
	A graph $H'$ is a \emph{restriction of $H$} on $G'$ 
	if $H'= H[\{v \times L(v): v \in V(G')\}.$ 
	Assume $G'$ has an $(H',L')$-coloring with an independent set $I'$ in $H'$ such that  
	$|I'|= |V(G)|-|V(F)|.$\\
	\indent A \emph{residual list assignment} $L^*$ of $F$ is defined by    
	$$L^*(x)=L(x)-\bigcup_{ux\in E(G)}\{c'\in L(x) : 
	(u,c)(x,c')\in M_{L,ux} \text{ and }(u,c)\in I'\}$$
	for each $x \in V(F).$ \\
	\indent A \emph{residual cover $H^*$} is defined by $H^*=H[\{x \times L^*(x): x \in V(F)\}].$
\end{defn}

From above definitions, we have the following fact. 

\begin{lem}\label{extend} 
	A residual cover $H^*$ is a cover of $F$ with an assignment $L^*. $
	Furthermore, if $F$ is $(H^*,L^*)$-colorable, then $G$ is $(H,L)$-colorable. 	
\end{lem}
\begin{pr}
	One can check from the definition of a cover and residual cover that 
	$H^*$ is a cover of $F$ with an assignment $L^*.$\\ 
	Suppose that  $F$ is $H^*(L^*)$-colorable. 
	Then $H^*$ has an independent set $I^*$ with $|I^*|= |F|.$ 
	It follows from Definition \ref{residual} that no edges connect $H^*$ and $I'.$ 
	Additionally, $I'$ and $I^*$ are disjoint. 
	Altogether, we have that $I=I'\cup I^*$ is an independent set in $H$   
	with $|I|=(|V(G)|-|V(F))+|V(F)|=|V(G)|.$ 
	Thus $G$ is $H(L)$-colorable. 	
\end{pr}

From now on, let  $G$ be a minimal non DP-$4$-colorable graph. 

\begin{lem}\label{lem1} 
	Every vertex in $G$ has degree at least $4.$
\end{lem}
\begin{pr} 
	Suppose to the contrary that $G$ has a vertex $x$ degree at most $3.$ 
	Let $L$ be a $4$-assignment 
	and let $H$ be a cover of $G$ such that $G$ has no $(H,L)$-coloring.       
	By the minimality of $G,$ the subgraph $G'=G-x$ admits 
	where $L'$ ($H'$) is a restriction of $L$ ($H$) in $G'.$ 
	Thus there is an independent set $I'$ with $|I'|= |G'|$ in $H'.$  
	Consider a residual list assignment $L^*$ on $x.$ 
	Since $|L(x)| = 4$ and $d(x)\leq 3$, we obtain $|L^*(x)|\geq 1.$  
	Clearly, $\{(x,c)\}$ is an independent set in $G[x]$ where $c \in L^*(x).$  
	Thus $G[x]$ is $(H^*,L^*)$-colorable. 
	It follows from Lemma \ref{extend} that the graph $G$ is $(H,L)$-colorable, a contradiction. 
\end{pr}

\begin{lem}\label{lem2} 
	If $F$ is an induced subgraph of $G$ obtained from 
	a cycle $x_1x_2\ldots x_m x_1$ and $k$ chords 
	$x_1x_{i_1}, x_1x_{i_2},\ldots,x_1x_{i_k},$ 
	then $d(x_1) \geq k+4$ or  $d(x_i) \geq 5$ for some $i \in \{2,3,\ldots,m.\}$
\end{lem}

\begin{pr} Suppose to the contrary that $d(x_1) \leq k+3$ and  $d(x_i) \leq 4$ for $i =2,3,\ldots,m.$ 
	Let $L$ be a $4$-assignment 
	and let $H$ be a cover of $G$ such that $G$ has no $(H,L)$-coloring.   
	By the minimality of $G,$ 
	the subgraph $G'=G-F$ admits an $(H',L')$-coloring 
	where $L'$ ($H'$) is a restriction of $L$ ($H$) in $G'.$ 
	Thus there is an independent set $I'$ with $|I'|= |G'|$ in $H'.$  
	Consider a residual list assignment $L^*$ on $x.$ 
	Since $|L(v)|= 4$ for every $v\in V(G),$ 
	we have $|L^*(x_1)|,$ $|L^*(x_{i_1})|,$ $|L^*(i_{k_2})|,$ $\ldots,$ $|L^*(x_{k_k})| \geq 3$ 
	and $|L^*(x_i)| \geq 2$ for remaining vertices in $H.$ 
	Let $H^*$ be an residual cover of $F.$ 
	We can  choose a color $c$ from  $L^*(x_1)$ such that 
	$|L^*(x_{m})-\{c':(x_1,c)(x_{m},c')\in H^*\}|\geq 2.$ 
	By choosing colors of $x_2$, $x_3,\dots, x_{m}$ in this order, 
	we obtain an independent set $I^*$ with $|I^*|=m=|F|.$ 
	Thus $F$ is $(H^*,L^*)$-colorable. 
	It follows from Lemma \ref{extend} that the graph $G$ is $(H,L)$-colorable, a contradiction. 
\end{pr}

By Lemma \ref{lem2}, we obtain the lower bound of the number of incident faces of a $6$-vertex that are incident to at least two $5^+$-vertices.  
\begin{cor}\label{cor4}  
	Every $6$-vertex $v$ in $G$  has at least two incident faces of $v$ that are incident to at least two $5^+$-vertices. 
\end{cor}

\section{Main Result}
\begin{thm} \label{thm1} Every planar graph without $i$-cycles is adjacent simultaneously to $j$-cycles and $k$-cycles is DP-$4$-colorable when $\{i,j,k\}=\{3,4,5\}$.
\end{thm}
\begin{pr}
\indent Suppose that $G$ is a minimal counterexample. Then each vertex in $G$ is a $4^+$-vertex by Lemma \ref{lem1}. The discharging process is as follows. Let the initial charge of a vertex $u$ in $G$ be $\mu(u)=2d(u)-6$ and the initial charge of a face $f$ in $G$ be $\mu(f)=d(f)-6$. Then by Euler's formula $|V(G)|-|E(G)|+|F(G)|=2$ and by the Handshaking lemma, we have
$$\displaystyle\sum_{u\in V(G)}\mu(u)+\displaystyle\sum_{f\in F(G)}\mu(f)=-12.$$
\indent Now, we establish a new charge $\mu^*(x)$ 
for all $x\in V(G)\cup F(G)$ by transferring charge from one element to another
and the summation of new charge $\mu^*(x)$ remains $-12$. 
If the final charge  $\mu^*(x)\geq 0$ for all $x\in V(G)\cup F(G)$, 
then we get a contradiction and the proof is completed.\\
\indent Let $w(v \rightarrow f)$ be the charge transferred from 
a vertex $v$ to an incident face $f.$ 
A $4$-vertex $v$ is \emph{flaw} if $v$  is a $(3,3,5,5^+)$-vertex. 
The discharging rules are as follows.\\
\textbf{(R1)} Let $f$ be a $3$-face.\\
\indent \textbf{(R1.1)} For a $4$-vertex $v$,\\
$
 w(v \rightarrow f) =
  \begin{cases}
     0.6,        & \text{if } v \text{ is flaw and } f \text{ is a } (4,5^+,5^+)\text{-face},\\
   0.8,        & \text{if } v \text{ is flaw and } f \text{ is a } (4,4,5^+)\text{-face},\\
   1,        & \text{otherwise.}
  \end{cases}
$\\
\indent \textbf{(R1.2)} For a $5^+$-vertex $v$,\\ 
$
 w(v \rightarrow f) =
  \begin{cases}
   1.4,        & \text{if } f  \text{ is a } (4,4,5^+)\text{-face with two incident flaw vertices},\\
     1.2,       & \text{if } f  \text{ is a } (4,4^+,5^+)\text{-face with exactly one flaw vertex},\\
   1,        & \text{otherwise.}
  \end{cases}
$\\
\textbf{(R2)} Let $f$ be a $4$-face.\\
\indent  For a $4^+$-vertex $v$,
$w(v \rightarrow f) = 0.5$.\\
\textbf{(R3)} Let $f$ be a $5$-face.\\
\indent \textbf{(R3.1)} For a $4$-vertex $v$,\\ 
$
 w(v \rightarrow f) =
  \begin{cases}
     0,        & \text{if } v  \text{ is  a flaw vertex with four } 4\text{-neighbors},\\
     0.1,        & \text{if } v  \text{ is  a flaw vertex with exactly one } 5^+\text{-neighbor},\\
   0.2,        & \text{if } v  \text{ is  a flaw vertex with at least two } 5^+\text{-neighbors},\\
   0.2,        & \text{if } v  \text{ is a } (3,5,5,4)\text{-vertex},\\
       1/3,        & \text{otherwise}.
  \end{cases}
$\\
\indent \textbf{(R3.2)} For a $5$-vertex $v$,\\
$
 w(v \rightarrow f) =
  \begin{cases}
   0.7,        & \text{if } f \text{ is a } (4,4,4,4,5)\text{-face with five adjacent }   4^-\text{-faces},\\
   0.6,        & \text{if } f \text{ is a } (4,4,4,4,5)\text{-face and at least one adjacent } 5^+\text{-face},\\
   0.4,        & \text{if } f \text{ is a } (4,4,4,5,5^+)\text{-face with both incident } 5^+\text{-vertices are adjacent}\\ 
   0.3,        & \text{otherwise.}
  \end{cases}
$\\
\indent \textbf{(R3.3)} For a $6^+$-vertex $v$,\\ 
$
 w(v \rightarrow f) =
  \begin{cases}
   0.8,        & \text{if } f \text{ is a } (4,4,4,4,6^+)\text{-face},\\
       0.4,        &\text{if } f \text{ is incident to a }  5^+\text{-vertex other than } v.
  \end{cases}
$\\
It remains to show that resulting $\mu^*(x)\geq 0$ for all $x\in V(G)\cup F(G)$. 
Moreover, it is evident for each $6^+$-face $f$ that $\mu^*(f)\geq 0$.\\

\indent \textit{CASE 1:} A $4$-vertex $v.$\\
\indent We use  (R1.1), (R2), and (R3.1) to prove this case.\\ 
\indent \textit{SUBCASE 1.1:} $v$ is a flaw vertex, that is a $(3,3,5,5^+)$-vertex.\\
\indent If each adjacent vertex of $v$ is a $4$-vertex, then  $\mu^*(v) \geq \mu(v) - 2 \times1= 0$. 
If $v$ is adjacent to exactly one $5^+$-vertex, then we obtain $\mu^*(v) \geq \mu(v) - 1-0.82\times0.1 = 0$. 
Now, assume $v$ is adjacent to at least two $5^+$-vertices. 
If $v$ is incident to a $(4,5^+,5^+)$-face, then $\mu^*(v) \geq \mu(v) - 1-0.6-2\times0.2= 0,$ 
otherwise,  $\mu^*(v) \geq \mu(v) - (2\times0.8+2\times0.2) = 0.$\\
\indent \textit{SUBCASE 1.2:} $v$ is not a flaw vertex.\\
\indent If $v$ is not incident to any $3$-face, then  $\mu^*(v)\geq \mu(v) - (4\times0.5 = 0.$ 
It follows from Proposition \ref{prop1} that  
$v$ with an incident $3$-face has at most one incident $4$-face.   
If $v$ is a $(3,5,5,4)$-vertex, then  $\mu^*(v)\geq \mu(v) - 1 -0.5-2\times0.2> 0.$ 
$v$ is a $(3,6^+,5^+,4)$-vertex, then $\mu^*(v)\geq \mu(v) - 1 -0.5-\frac{1}{3}> 0.$ 
If $v$ is incident to exactly one $3$-face but is not incident to any $4$-face, 
then $\mu^*(v)\geq \mu(v) - 1 - 3\times 1/3 = 0.$ 
If $v$ is incident to two $3$-faces, then the remaining two incident faces are $6^+$ by Proposition \ref{prop1}. 
Thus $\mu^*(v) \geq \mu(v) - 2\times 1= 0.$\\
\indent\textit{CASE 2:} Consider a $5$-vertex $v.$\\ 
\indent We use  (R1.2), (R2.2), and (R3.2) to prove this case.\\
\indent If $v$ is incident to some $3$-face in $C(3,3)$, then two  incident faces of $v$ are $6^+$-faces by Proposition \ref{prop1}. 
Thus  $\mu^*(v) \geq \mu(v) - 1.4-2\times1 > 0.$ 
Now, assume $v$ is not incident to $C(3,3)$. This implies that $v$ is adjacent to at most two $3$-faces.\\
\indent \textit{SUBCASE 2.1:} $v$ is  incident to at most one $3$-face. \\
\indent If $v$ is not incident to any $4$-face, then each incident $5$-face of $v$ is adjacent to some $5^+$-face. 
Thus $\mu^*(v) \geq \mu(v) - 1.4-4\times0.6 > 0.$  
If $v$ is incident to at least one $4$-face, then   $\mu^*(v) \geq \mu(v) - 1.4-3\times0.7-0.5 = 0.$\\
\indent \textit{SUBCASE 2.2:}  $v$ is incident to two $3$-faces. \\
Let $v$ be incident to faces $f_1,f_2,$ $\ldots,$ $f_5$ in a cyclic order 
where $f_1$ and $f_3$ are $3$-faces. 
It follows from Proposition \ref{prop1} that $f_2$, $f_4$, and $f_5$ are $5^+$-faces. 
Assume $f_1$ and $f_3$ are $(4,4,5)$-faces. 
It follows from Corollary \ref{cor4} that $f_2$ is incident to at least two non-adjacent $5^+$-vertices.  
If either $f_4$ or $f_5$ is incident to one $5^+$-vertex, 
then the other is incident to at least non-adjacent two $5^+$-vertices by Corollary \ref{cor4}. 
Thus  $\mu^*(v) \geq \mu(v) - 2\times1.4-0.6-2\times0.3 > 0.$ 
If both $f_4$ and $f_5$ are incident to at least two $5^+$-vertices, then $\mu^*(v) \geq \mu(v) - 2\times1.4-2\times0.4-0.3 > 0.$ 
If  $f_1$ (or $f_3$) is a $(5,5^+,5^+)$-face, then $f_2$ and $f_5$ (or $f_4$) are incident to at least two $5^+$-vertices. \\
Thus $\mu^*(v) \geq \mu(v) - 1.4-1-0.6-2\times0.4 > 0.$ \\
\indent Assume that both $f_1$ and $f_3$ are $(4,5,5^+)$-faces where $5^+$-vertex other than $v$ are $x$ and $y$. 
If $x$ and $y$ are not in $f_2$, then both $f_4$ and $f_5$ are adjacent to at least two $5^+$-vertices. 
Thus $\mu^*(v) \geq \mu(v) - 2\times1.2-0.7-2\times0.4 > 0.$
If $x$ and $y$ are  in $f_2$, then  $\mu^*(v) \geq \mu(v) - 2\times1.2-2\times0.6-0.3 > 0.$  
If exactly one of $x$ and $y$ is in $f_2,$ then $\mu^*(v) \geq \mu(v) - 2\times1.2-0.6-2\times0.4 > 0.$\\
\indent Assume that exactly one of $f_1$ and $f_3$ is a $(4,5,5^+)$-face where $5^+$-vertex other than $v$ is $w$ 
and $w$ is in $f_2.$   
It follows from Corollary \ref{cor4} that  at least one of  $f_4$ and $f_5$ is incident to at least two  $5^+$-vertices. 
If $f_4$ or $f_5$ is incident to at least two $5^+$-vertices and two of them are not adjacent, 
then $\mu^*(v) \geq \mu(v) - 1.4-1.2-0.6-0.4-0.3 > 0,$ 
otherwise $f_4$ and $f_5$ are incident to at least two $5^+$-which implies $\mu^*(v) \geq \mu(v) - 1.4-1.2-3\times0.4 > 0.$ \\
\indent Assume that exactly one of $f_1$ and $f_3$ is a $(4,5,5^+)$-face where $5^+$-vertex other than $v$ is $w$ 
and $w$ is in $f_5.$   
It follows that one of $f_2$ and $f_4$ is incident to  at least two $5^+$-vertices. 
If $f_2$ is incident at least two $5^+$-vertices, then two of them are not adjacent. Thus $\mu^*(v) \geq \mu(v) - 1.4-1.2-0.6-0.4-0.3 >0.$ 
If  $f_5$ is  incident at least two $5^+$-vertices, then two of them are not adjacent or two of them are not not adjacent. 
Thus $\mu^*(v) \geq \mu(v) - 1.4-1.2-0.7-0.4-0.3 =0.$ \\
\indent \textit{CASE 3:} A $6$-vertex $v.$\\ 
\indent We use  (R1.2), (R2.2), and (R3.3) to prove this case.\\
\indent \textit{SUBCASE 3.1:} $v$ is incident to $C(3,3)$. \\
\indent  It follows from Proposition \ref{prop1}(b) that $v$ is incident to at least two $6^+$-faces. 
Then $\mu^*(v) \geq \mu(v) - 4\times1.4>0.$\\
\indent \textit{SUBCASE 3.2:} $v$ is not incident to  $C(3,3)$.\\ 
\indent It follows that $v$ is incident at most three $3$-faces.
  If $v$ is incident to at most two $3$-faces, then $\mu^*(v) \geq \mu(v) - 2\times1.4-4\times0.8= 0.$ 
  Assume $v$ is incident to three $3$-faces. 
By Proposition \ref{prop1}(a), $v$ is not  incident to any $4$-face. 
If $v$ is incident to at least one $6^+$-face, then $\mu^*(v) \geq \mu(v) - 3\times1.4-2\times0.8>0.$ 
Now, it remains to consider the case that  $v$ is a $(3,3,3,5,5,5)$-vertex. 
By Corollary \ref{cor4},  two of incident faces of $v$, say $f_1$ and $f_2$ (not necessary adjacent),  
are incident  to at least two $5^+$-vertices. 
If $f_1$ and $f_2$ are $5^+$-faces, then $\mu^*(v) \geq \mu(v) - 3\times1.4-0.8-2\times0.4>0.$ 
If $f_1$ or $f_2$ is a $3$-face, then $\mu^*(v) \geq \mu(v) - 2\times1.4-1.2-2\times0.8-0.4=0.$  \\
\indent \textit{CASE 4:} A $k$-vertex $v$ with $k\geq 7.$\\
\indent We use  (R1.2), (R2.2), and (R3.3) to prove this case.\\
\indent \textit{SUBCASE 4.1:} $v$ is incident to $C(3,3)$. \\
\indent It follows from Proposition \ref{prop1} that  $v$ is incident to at least two $6^+$-faces. 
Thus  $\mu^*(v) \geq \mu(v) - (k-2)\times1.4>0$ for $k\geq 6$.\\
\indent \textit{SUBCASE 4.2:} $v$ is not incident to  $C(3,3)$.\\
\indent It follows that $v$ is incident to at most $\frac{d(v)}{2}$ $3$-faces. 
Thus  $\mu^*(v) \geq \mu(v)-\frac{k}{2}\times1.4-\frac{k}{2}\times0.8$ for $k\geq7$.\\
\indent \textit{CASE 5:}  A $3$-face.\\ 
\indent We use (R1.1) and (R1.2) to prove this case.\\   
\indent If each vertex of $f$ is not a flaw vertex, then $\mu^*(f) = \mu(f) +3\times1=0$.\\ 
 \indent Now, assume that  at least one incident $4$-vertex of $f$ is flaw.  
If $f$ is a $(4,4,5^+)$-face with two incident flaw vertices, then $\mu^*(f) = \mu(f) +1.4+2\times0.8=0.$ 
If $f$ is a $(4,4,5^+)$-face with exactly one incident flaw vertex, then $\mu^*(f) =\mu(f) +1.2+1+0.8=0.$  
If $f$ is a $(4,5^+,5^+)$-face, then $\mu^*(f) = \mu(f) +2\times1.2+0.6=0$.\\
\indent \textit{CASE 6:} A $4$-face $f.$\\ 
\indent We obtain  $\mu^*(f) \geq \mu(f)  +4\times0.5= 0$ by (R2). \\
\indent \textit{CASE 7:} A $5$-face $f.$\\
 \indent We use (R3.1), (R3.2), and (R3.3)  to prove this case.\\ 
\indent \textit{SUBCASE 7.1:} $f$ is incident to at least three $5^+$-vertices.\\
 \indent It follows that each $4$-vertex in $f$ is adjacent to at least one $5^+$-vertex. Then $\mu^*(f)\geq\mu(f)+ 3\times0.3+2\times0.1> 0.$ \\  
\indent \textit{SUBCASE 7.2:} $f$ is incident to two $5^+$-vertices, say $x$ and $y$.\\
\indent 
If $x$ and $y$ is not adjacent, then each $4$-vertex in $f$ is adjacent to at least one $5^+$-vertex. 
Additionally, one of them is adjacent to at least two $5^+$-vertices. Thus $\mu^*(f)\geq\mu(f)+ 2\times0.3+0.2+2\times0.1= 0.$ 
If $x$ and $y$ are adjacent, then two incident $4$-vertices of $f$ are adjacent to at least one $5^+$-vertex. 
Thus $\mu^*(f)\geq\mu(f)+ 2\times0.4+2\times0.1= 0.$\\  
\indent \textit{SUBCASE 7.3:} $f$ is incident to one $5^+$-vertex.\\
\indent If  $f$ is a $(4,4,4,4,6^+)$-face,  then two $4$-vertices in $f$ are adjacent to at least one $6^+$-vertex. 
Thus  $\mu^*(f)\geq\mu(f)+ 0.8+2\times0.1= 0.$ 
Assume $f$ is a $(4,4,4,4,5)$-face with five incident vertices $x_1,x_2,\ldots,x_5$ in a cyclic order with $d(x_5)=5.$ 
If  each adjacent face of $f$ is a $4^-$-face and both $x_1$ and  $x_4$ are flaw, 
then there is a a $(3,5,5^+)$-face adjacent to $f$ and incident to $x_5$  by Corollary \ref{cor4}. 
Consequently, $x_1$ or $x_4$ s adjacent to at least two $5^+$-vertices. 
Thus $\mu^*(f)\geq\mu(f)+ 0.7+0.2+0.1= 0.$  
If  each adjacent face of $f$ is a $4^-$-face but $x_1$ or $x_4$ is not flaw, 
then $\mu^*(f)\geq\mu(f)+ 0.7+0.2+0.1= 0.$ 
If one adjacent face of $f$ is a $5^+$-face, then one of $4$-vertex in $f$ is not either flaw or $(3,5,5,4)$-vertex. 
Thus $\mu^*(f)\geq\mu(f)+ 0.6+\frac{1}{3}+0.1> 0.$ \\
\indent \textit{SUBCASE 7.4:} $f$ is a $(4,4,4,4,4)$-face.\\
\indent  It follows from Lemma \ref{lem2} that  each adjacent $3$-face of $f$ is a $(4,4,5^+)$-face. 
This implies that each flaw vertex in $f$ is adjacent to at least two $5^+$-vertex. 
Consequently, each vertex in $f$ send charge at least $0.2$ to $f$. 
Thus  $\mu^*(f)\geq\mu(f)+ 5\times0.2= 0.$\\  
\indent This completes the proof. 
\end{pr}
\section{Acknowledgments}
\indent The first author is supported by Development and Promotion of Science and Technology
talents project (DPST).

\end{document}